\documentclass[12pt]{article}
\usepackage{amsfonts}
\usepackage{amssymb}
\usepackage{enumerate}
\usepackage[T1]{fontenc}

\newtheorem{ttt}{Theorem}[section]
\newtheorem{llll}[ttt]{Lemma}
\newtheorem{ccc}[ttt]{Claim}
\newtheorem{eee}[ttt]{Example}
\newtheorem{sss}[ttt]{Statement}
\newtheorem{ddd}[ttt]{Definition}
\newtheorem{qqq}[ttt]{Question}
\newtheorem{cccc}[ttt]{Corollary}

\newcommand{\bt}{\begin{ttt}}
\newcommand{\bl}{\begin{llll}}
\newcommand{\bc}{\begin{ccc}}
\newcommand{\bex}{\begin{eee}}
\newcommand{\bs}{\begin{sss}}
\newcommand{\bd}{\begin{ddd} \upshape}
\newcommand{\bq}{\begin{qqq}}
\newcommand{\bcor}{\begin{cccc}}

\newcommand{\bp}{\noindent\textbf{Proof. }}

\newcommand{\et}{\end{ttt}}
\newcommand{\el}{\end{llll}}
\newcommand{\ec}{\end{ccc}}
\newcommand{\eex}{\end{eee}}
\newcommand{\es}{\end{sss}}
\newcommand{\ed}{\end{ddd}}
\newcommand{\eq}{\end{qqq}}
\newcommand{\ecor}{\end{cccc}}

\newcommand{\ep}{\hspace{\stretch{1}}$\square$\medskip}

\newcommand{\lab}[1]{\label{#1}}

\newcommand{\RR}{\mathbb{R}}

\newcommand{\e}{\varepsilon} 
\newcommand{\om}{\omega}

\newcommand{\iB}{\mathcal{B}}
\newcommand{\iE}{\mathcal{E}}

\newcommand{\iI}{\mathcal{I}}

\newcommand{\iF}{\mathcal{F}}

\newcommand{\iN}{\mathcal{N}}

\begin{document}

\title{Less than $2^{\omega}$ many translates of a compact nullset may cover the real line}

\author{M\'arton Elekes\thanks{Partially supported by Hungarian Scientific
Foundation grant no. 37758 and F 43620. The first author is also grateful
to Fields Institute and E\"otv\"os J\'ozsef Scholarship.} \ and Juris
Stepr$\bar{\textrm{a}}$ns}

\maketitle 

\begin{abstract}
We answer a question of Darji and Keleti by proving that there exists
a compact set $C_0\subset\RR$ of measure zero such that for every perfect set
$P\subset\RR$ there exists $x\in\RR$ such that $(C_0+x)\cap P$ is
uncountable. Using this $C_0$ we answer a question of Gruenhage by showing
that it is consistent with $ZFC$ (as it follows e.g.~from
$\textrm{cof}(\iN)<2^\omega$) that less
than $2^\omega$ many translates of a compact set of measure zero can cover $\RR$.
\end{abstract}

\insert\footins{\footnotesize{MSC codes: Primary 28E15; Secondary 03E17,
03E35}}
\insert\footins{\footnotesize{Key Words: compact, measure zero, translate,
cover, continuum, perfect, consistent}}

\bigskip
\bigskip

\section*{Introduction}


\bigskip

The behaviour of the classical cardinal invariants in the Cicho\'n diagram is
very well described. See e.g.~the monograph \cite{BJ} for definitions and
details. We also follow the terminology of \cite{BJ} throughout the paper. 

The invariant we are interested in is \emph{$\textrm{cov}(\iN)$}, that is
the
least cardinal $\kappa$ for which it is possible to cover $\RR$ by $\kappa$
many nullsets (sets of Lebesgue measure zero), and also some variants of
$\textrm{cov}(\iN)$. There are two natural ways to modify this
definition. (See \cite{BJ} Chapter 2.6 and 2.7.) First,
\emph{$\textrm{cov}^*(\iN)$} is the least cardinal $\kappa$ for which it is
possible to cover $\RR$ by $\kappa$ many translates of some nullset. In other
words, $\textrm{cov}^*(\iN) = \min\{|A| \ | \ A\subset\RR, \exists N \in \iN, A+N=\RR
\}$, where $|A|$ is the cardinality of $A$, $\iN$ is the $\sigma$-ideal of
nullsets and $A+N=\{a+n\  |\  a\in A, n\in N\}$. The other possible modification
is \emph{$\textrm{cov}(c\iN)$}, that is the least cardinal $\kappa$ for which it is
possible to cover $\RR$ by $\kappa$ many compact nullsets. (At this point we
depart from the terminology of \cite{BJ} as this notion is denoted by
$\textrm{cov}(\iE)$ there.) It can be found in these two chapters of this
monograph that both $\textrm{cov}^*(\iN)<2^\omega$ and
$\textrm{cov}(c\iN)<2^\omega$ are consistent with $ZFC$.

G. Gruenhage posed the natural question whether
\emph{$\textrm{cov}^*(c\iN)<2^\omega$} is also consistent, that is, whether we
can consistently cover $\RR$ by less than continuum many translates of a
compact nullset.

The main goal of this paper is to answer this question in the affirmative via
an answer (in
$ZFC$) to a question of U. B. Darji and T. Keleti that is also interesting in
its own right.



We remark here that under $CH$ (the Continuum Hypothesis, or more generally
under $\textrm{cov}(\iN)=2^\omega$) the real line
obviously cannot be covered by less than $2^\omega$ many nullsets. Therefore
it is consistent that the type of covering we are looking for does not
exist.

So the interesting case is when the consistent inequality
$\textrm{cov}^*(\iN)<2^\omega$ holds. The nullset in this statement can obviously be
chosen to be $G_\delta$. So the content of Gruenhage's question
actually is whether this can be an $F_\sigma$ or closed or
compact nullset. We formulate the strongest version.

\bq 
(Gruenhage) Is it consistent that there exists a compact set $C\subset
\RR$ of Lebesgue measure zero and $A\subset \RR$ of cardinality less than
$2^\omega$ such that $C+A=\RR$?
\eq

For example Gruenhage showed that no such covering is possible if $C$ is the
usual ternary Cantor set (see \cite{DK} and for another motivation of this
question see \cite{GL}).

Working on this question Darji and Keleti \cite{DK} introduced
the following notion.

\bd 
(Darji - Keleti) Let $C\subset\RR$ be arbitrary. A set $P\subset\RR$ is
called a \emph{witness for} $C$ if $P$ is perfect and for every translate
$C+x$ of $C$ we have that $(C+x)\cap P$ is countable.
\ed

Obviously, if there is a witness $P$ for $C$ then less than $2^\om$ many
translates of $C$ cannot cover $P$, so they cannot cover $\RR$. Motivated by a
question of D. Mauldin, who asked what can be said if $C$ is of
Hausdorff dimension strictly less than 1, Darji and Keleti proved the
following. 

\bt
(Darji - Keleti) If $C\subset\RR$ is a compact set of packing dimension 
$\dim_p(C)<1$ then there is a witness for $C$, and
consequently less than $2^\om$ translates of $C$ cannot cover $\RR$.
\et

They posed the following question, an affirmative answer to which would also
answer the original question of Gruenhage in the negative.

\bq
(Darji - Keleti) Is there a witness for every compact set $C\subset\RR$ of
Lebesgue measure zero?
\eq

We will answer this question in the negative, which still leaves the original
question of Gruenhage open. Then we will show that using the same ideas it is
also possible to give an affirmative answer to Gruenhage's question.

\section{Answer to the question of Darji and Keleti}

The following set is fairly well known in geometric measure theory, as it is
probably the most natural example of a compact set of measure zero but of
Hausdorff and packing dimension 1. It was investigated for example by Erd\H os
and Kakutani \cite{EK}.

\bd 
Denote
\[
C_0 = \left\{ \left. \sum_{n=2}^{\infty} \frac{d_n}{n!} \ \right| \ d_n
\in\{0,1,\dots, n-2\}\ \forall n \right\}.
\]
\ed

Think of $d_n$ as digits with ``increasing base''; then all but countably
many $x\in[0,1]$ have a unique expansion
\[
x=\sum_{n=2}^\infty \frac{x_n}{n!},
\]
where $x_n\in \{0,1,\dots, n-1\}$ for every $n=2,3,\dots$

The set of real numbers with prescribed first $n$ digits is a closed interval.
Let us call these \emph{basic intervals of level n} and denote this collection by
$\iB_n$. Let $\iB=\cup_{n=0}^\infty \iB_n$ be the set of all \emph{basic intervals}.
So $\iB$ forms a tree under inclusion, the $n^{th}$ level of
which is $\iB_n$, which consists of $n!$ nonoverlapping intervals of length
$\frac{1}{n!}$.

Using the above expansion it is easy to see that $C_0\subset\RR$ is a compact
set of Lebesgue measure zero.

\bigskip

The following theorem answers the question of Darji and Keleti.

\bt\lab{intersects}
For every perfect set $P\subset\RR$ there exists a translate $C_0+x$ of the
compact nullset $C_0$ such that $(C_0+x)\cap P$ is uncountable.
\et

\bp
We will show that there exists a perfect set $Q\subset P$ and $y\in \RR$ such
that $Q+y\subset C_0$. This is clearly sufficient, as then for $x=-y$ we get
$Q\subset (C_0+x)\cap P$, so this intersection is uncountable. First we will
construct $Q$ via a dyadic tree of basic intervals, then we will construct
the ``digits'' of $y$.

By translating $P$ if necessary we can assume that $P$ intersects
$(0,\frac{1}{5!})$.  Instead of $P$ we may as well work with any perfect
subset of it, for if we find the set $Q$ inside this subset then this $Q$ also
works for $P$. So we may find a perfect subset of $P$ in $[0,\frac{1}{5!}]$
and therefore we can assume that
$P$ itself is inside $[0,\frac{1}{5!}]$.  Moreover, as $P$ is uncountable and the endpoints
of the basic intervals form a countable set, we can find a perfect subset of
$P$ that is disjoint from the set of endpoints (we used here twice the well known fact that
every uncountable Borel set contains a perfect set). Therefore we can assume
that $P$ itself is disjoint from the endpoints.

Now we recursively pick an increasing sequence of levels $(l_k)_{k=0}^\infty$
and for every $k$ choose a set $\iI_{l_k}\subset\iB_{l_k}$ of size $2^k$ such
that
\begin{enumerate}
\item for each $I\in \iI_{l_k}$ there are exactly two intervals in
$\iI_{l_{k+1}}$ (at level $l_{k+1}$) that are contained in $I$, the so called
\emph{successors} of $I$
\item\lab{intersect} for each $I\in \iI_{l_k}$ we have $I\cap P \neq
\emptyset$
\item\lab{far} $l_k \geq 2^{k+2}+1$.
\end{enumerate}

The recursion is carried out as follows. Fix $p_0\in P$. Let $l_0=5$, and
at level 5 we pick the (unique) basic interval $I$ containing $p_0$ in its
interior. Let $\iI_{l_0}=\iI_5=\{I\}$. The recursion step is as follows. As
$P$ is disjoint from the endpoints of the basic intervals, each interval
$I\in\iI_{l_k}$ (at level $l_k$) contains some point $p_I\in P$ in its
interior by condition \ref{intersect}. As $P$ is perfect, we can choose a
distinct point $p'_I\in P$ in $I$. We can find a large enough $n$ such that
the $2^{k+1}$ distinct points $p_I$ and $p'_I \ (I\in\iI_{l_k})$ are all
separated by $\iB_n$. Define
\[
l_{k+1}=\max\{n,2^{k+3}+1\}.
\]
Clearly, condition \ref{far} is also satisfied.

Let $\iI_{l_{k+1}}$ be the subcollection of $\iB_{l_{k+1}}$ consisting of the
$2^{k+1}$ basic intervals containing all the points $p_I$ and $p'_I\
(I\in\iI_{l_k})$. This recursion clearly provides a system of intervals
satisfying the required properties.

Now we can define
\[
Q=\bigcap_{k=0}^\infty \bigcup \iI_{l_k}.
\]

Let us extend this tree of intervals to the intermediate levels in the
natural way, that is, for every $I\in\iI_{l_k}$ and successor
$J\in\iI_{l_{k+1}}$ and every $n\in(l_k,l_{k+1})$ let us add to the tree the unique basic interval of level $n$ that is contained in $I$ and contains $J$. For $n=2,3,4,5$
define $\iI_n=\{[0,\frac{1}{n!}]\}$. Hence we get $\iI_n$ for every
$n=2,3,\dots$ so that
\[
\bigcap_{n=2}^\infty \bigcup \iI_n=\bigcap_{k=0}^\infty \bigcup \iI_{l_k}=Q.
\]

Our next goal is to define $y=\sum_{n=2}^\infty \frac{y_n}{n!}$ so that
$Q+y\subset C_0$. Define $y_2=y_3=y_4=y_5=0$. For every $n\geq 6$ there
exists $k$ such that $l_k<n\leq l_{k+1}$. Clearly, the size of $\iI_n$ is
$2^{k+1}$, and $Q\subset \cup \iI_n$. This means that there are at most
$2^{k+1}$ possible values for $q_n$, where $q\in Q$ and $q=\sum_{n=2}^\infty
\frac{q_n}{n!}$ (we do not have to worry about nonunique expansions, as
$Q\subset P$ so $Q$ is disjoint from the endpoints of the basic intervals).
For every such $q_n$ there are at most two values of $m$ such that
$q_n+m\in\{n-2,n-1\}$. Hence altogether there are at most $2\cdot2^{k+1}$ such
``bad'' values, so if $n-1>2\cdot2^{k+1}$ then we can fix a
$y_n\in\{0,1,\dots,n-2\}$ such
that $q_n+y_n\notin\{n-2,n-1\}$ for every possible $q_n$. But our requirement
on $n$ and $k$, namely $n-1>2\cdot2^{k+1}$,
is clearly satisfied as $n>l_k\geq 2^{k+2}+1$ by condition \ref{far}.

So we can define
\[
y=\sum_{n=2}^\infty \frac{y_n}{n!}
\]
so that for every $n\geq 6$ we have $y_n\in\{0,1,\dots,n-2\}$ and that for every
$q\in Q$ with $q=\sum_{n=2}^\infty \frac{q_n}{n!}$ we have
$q_n+y_n\notin\{n-2,n-1\}$. We claim that $Q+y\subset C_0$, which will
complete the proof. Fix $q\in Q$ with $q=\sum_{n=2}^\infty \frac{q_n}{n!}$,
then
\[
q+y = \sum_{n=2}^\infty \frac{q_n+y_n}{n!} = \sum_{n=2}^\infty
\frac{q_n+y_n-n\varepsilon_n}{n!} + \sum_{n=2}^\infty
\frac{n\varepsilon_n}{n!},
\]
where $\e_n$ is the ``carried digit'', so
\[
\e_n=\left\{ \begin{array}{ll} 0 & \textrm{if } q_n+y_n \leq n-1\\ 
                               1 & \textrm{otherwise}.
             \end{array}  \right.
\]
Continuing the above calculation we get
\[
q+y = \sum_{n=2}^\infty \frac{q_n+y_n-n\varepsilon_n}{n!} + \sum_{n=2}^\infty
\frac{\varepsilon_n}{(n-1)!} = \sum_{n=2}^\infty
\frac{q_n+y_n-n\varepsilon_n}{n!} + \sum_{n=1}^\infty
\frac{\varepsilon_{n+1}}{n!} =  
\]
\[
= \e_2 + \sum_{n=2}^\infty\frac{q_n+y_n-n\varepsilon_n+\e_{n+1}}{n!} =
\sum_{n=2}^\infty \frac{q_n+y_n-n\varepsilon_n+\e_{n+1}}{n!},
\]
since $\e_2=0$ by e.g. $y_2=0$. We now check that for every $n\geq 2$ the numerator
$q_n+y_n-n\varepsilon_n+\e_{n+1} \in \{0,1,\dots,n-2\}$, which shows that
$q+y\in C_0$. For $n<6$ this is clear, as $y_n=0$ and also $q_n=0$ by the
assumption $P\subset[0,\frac{1}{5!}]$. For
$n\geq 6$ recall that $q_n\leq n-1$, $y_n\leq n-2$, so $q_n+y_n \leq 2n-3$ and
also that $q_n+y_n\notin\{n-2,n-1\}$. We separate the cases $\e_n=0$ and
$\e_n=1$. If $\e_n=0$, then $q_n+y_n \leq n-1$, but then also $q_n+y_n \leq
n-3$. Therefore $q_n+y_n-n\varepsilon_n+\e_{n+1} =
q_n+y_n+\e_{n+1} \leq n-2$, and we are
done. On the other hand,
if $\e_n=1$, then
$q_n+y_n-n\varepsilon_n \leq n-3$, so $q_n+y_n-n\varepsilon_n+\e_{n+1} \leq
n-2$, so this case is also done. This completes the proof.
\ep

\section{Answer to the question of Gruenhage}

Now we answer the original question of Gruenhage. Recall that \emph{$\textrm{cof}(\iN)$} is the minimal cardinality of a family $\iF\subset\iN$ for which every nullset is contained in some member of $\iF$.

\bt\lab{less}
$\RR$ can be covered by $\textrm{cof}(\iN)$ many translates of $C_0$, consequently $\textrm{cov}^*(c\iN) \leq \textrm{cof}(\iN)$.
\et

\bp
It is clearly sufficient to cover the unit interval.

Fix
$f:\omega\setminus\{0,1\}\rightarrow\omega\setminus\{0\}$. A set of the form
$S=\Pi_{n=2}^\infty A_n$, where $A_n\subset \{0,1,\dots,n-1\}$ is of
cardinality at most $f(n)$ for every $n$ is called an
\emph{$f$-slalom}. Suppose $\lim_{n\to\infty}f(n) =\infty$, $f(2)=f(3)=f(4)=f(5)=1$ and
also that $f(n)<\frac{n-1}{2}$ for every $n\geq 6$. Combining \cite{BJ} Thm
2.3.9 and \cite{GL} Thm 2.10 we obtain a cover of $\Pi_{n=2}^\infty
\{0,1,\dots,n-1\}$ by $\textrm{cof}(\iN)$ many $f$-slaloms $\{S_\alpha \left| \
\alpha< \textrm{cof}(\iN) \}\right.$. (We actually obtain a cover of
$\omega^\omega$ first, which is a larger space, so it is trivial to restrict this cover to get a cover of $\Pi_{n=2}^\infty
\{0,1,\dots,n-1\}$. Moreover, \cite{BJ} works with
inclusion mod
finite, but that makes no difference, as we can replace each slalom by
countably many slaloms to get around this difficulty.) For a slalom $S$ define 
\[
S^* = \left\{ \left. \sum_{n=2}^\infty \frac{s_n}{n!} \right|\
\left(s_n\right)_{n=2}^\infty \in S \right\}.  
\]
Clearly $\{S^*_\alpha \left| \ \alpha< \textrm{cof}(\iN) \}\right.$ covers the
unit interval. 
The following lemma will complete the proof of the theorem.

\bl
Let $f$ be as above and $S$ be an $f$-slalom. Than there exists $y\in\RR$ such
that $S^*+y\subset C_0$.
\el

\bp
The proof is based on the ideas used in Theorem \ref{intersects}. $S^*$ plays
the role of $Q$. Our goal is to define
\[
y=\sum_{n=2}^\infty \frac{y_n}{n!}
\]
with $y_n\in\{0,1,\dots,n-2\}$ so that for every
$\left(s_n\right)_{n=2}^\infty \in S$ and for every $n \geq 6$ we have
$s_n+y_n\notin\{n-2,n-1\}$. But this is clearly possible by our assumptions on
$f$, as there are at most $f(n)<\frac{n-1}{2}$ possibilities for $s_n$, hence
there are two consecutive values excluded, and so we can find a suitable
$y_n\in\{0,1,\dots,n-2\}$.

Then by the same calculation as in the last part of the proof of Theorem
\ref{intersects} we check that
\[
\sum_{n=2}^\infty \frac{s_n+y_n}{n!} \in C_0.
\]
This completes the proof of the lemma.
\ep

Hence for every $\alpha < \textrm{cof}(\iN)$ there exists $y_\alpha$ such that
$S_\alpha^*+y_\alpha\subset C_0$, but than for $x_\alpha=-y_\alpha$ we have
$S_\alpha^*\subset C_0+x_\alpha$, so we obtain a cover of the unit interval by
$\textrm{cof}(\iN)$ many translates of $C_0$ and therefore the proof of
Theorem \ref{less} is also complete.
\ep

\bcor\lab{covers}
It is consistent that less than continuum many translates of a compact set of
measure zero cover the real line, that is, $\textrm{cov}^*(c\iN) < 2^\omega$
is consistent.
\ecor

\bp
$\textrm{cof}(\iN)$ is consistently less than the continuum \cite{BJ} p.~388.
\ep

\section{Remarks and open problems}

There are alternative ways to prove Theorem \ref{covers}. Using
Theorem \ref{intersects} one can directly use Sacks forcing to show
that the translates of $C_0$ by the ground model reals cover
$\RR$. Another method is to use the so called $CPA$ axiom (see
\cite{CP}). However, in all these cases $\textrm{cof}(\iN)$ is less
than the continuum. This is not necessary, as in the Laver model
$\textrm{cof}(\iN)=2^\omega=\omega_2$ but it can be shown that the
above slalom argument still works, that is, using the above $f$ the
space $\Pi_{n=2}^\infty \{0,1,\dots,n-1\}$ can be covered by
$\omega_1$-many $f$-slaloms in the Laver model. (This can be derived
from the so called \emph{Laver property}, see \cite{BJ}.) Therefore
$\textrm{cov}^*(c\iN)$ is not the same as $\textrm{cof}(\iN)$.





\bq
Is $\textrm{cov}^*(c\iN)$ equal to one of the known cardinal invariants?
\eq

Another natural question is the following. In most cases the values of the
classical cardinal invariants in the Cicho\'n diagram remain the same if we
replace $\RR$ with an arbitrary (uncountable) Polish space. However, in case
of $\textrm{cov}^*(c\iN)$ the situation is not clear. The authors were able to
reprove the results of this paper in case of the Cantor group or more
generally for countable products of finite discrete groups equipped with the
Haar measure, but not in the general case.

\bq
Can it be shown (without resorting
to extra set theoretic axioms) that there is an uncountable locally compact
Polish group G with Haar measure $\mu$ such that for every compact set $C\subset
G$ with $\mu(C)=0$ and every $A\subset G$ of cardinality less than $2^\omega$
we have $C+A\neq G$?
\eq


As for the cardinal invariants, even the following is open.

\bd
If $G$ is an uncountable locally compact Polish group with Haar
measure $\mu$ then \emph{$\textrm{cov}^*_G(c\iN)$} is the smallest cardinal
$\kappa$ for
which it is possible find a compact set of $\mu$-measure zero and cover $G$ by
$\kappa$ many translates of it.
\ed

\bq
Is it true that for every $G_1$ and $G_2$ uncountable locally compact Polish groups $\textrm{cov}^*_{G_1}(c\iN) = \textrm{cov}^*_{G_2}(c\iN)$?
\eq

\bigskip

\noindent
\textsc{R\'enyi Alfr\'ed Institute, Re\'altanoda u. 13-15. Budapest
1053, Hungary}

\textit{Email address}: \verb+emarci@renyi.hu+


\bigskip

\noindent
\textsc{Department of Mathematics, York University,
Toronto, Ontario M3J 1P3, Canada }

\textit{Email address}: \verb+steprans@mathstat.yorku.ca+

\end{document}